\documentclass{amsart}
\usepackage{graphicx}
\usepackage{algorithm}
\usepackage{algpseudocode}
\usepackage{booktabs}
\usepackage{multirow}
\usepackage{arydshln}
\usepackage[table,xcdraw]{xcolor}
\usepackage{tabu}
\usepackage{amssymb}
\usepackage{mathtools}
\usepackage{mathrsfs}
\usepackage{amsthm}
\numberwithin{equation}{section}
\numberwithin{figure}{section}
\numberwithin{table}{section}
\usepackage{bbm}
\usepackage{subfig}
\usepackage{enumerate}
\usepackage[section]{placeins}

\usepackage{ifpdf}
\ifpdf
\DeclareGraphicsExtensions{.pdf,.eps,.jpg,.png}	
\usepackage[suffix=]{epstopdf}
\fi
\usepackage{xcolor}
\usepackage[utf8]{inputenc}
\usepackage{hyperref}
\hypersetup{hidelinks}

\long\def\MSC#1\EndMSC{\def\arg{#1}\ifx\arg\empty\relax\else
	{\narrower\noindent%
		{2020 Mathematics Subject Classification}: #1\\} \fi}

\long\def\KEY#1\EndKEY{\def\arg{#1}\ifx\arg\empty\relax\else
	{\narrower\noindent%
		Keywords: #1\\}\fi}

\theoremstyle{plain}

\theoremstyle{definition}

\theoremstyle{remark}

\begin{document}
\title[A hybrid parametric modelling of ultrasound transmission]{A hybrid time-frequency parametric modelling of medical ultrasound signal transmission}




\author[C.~Razzetta]{Chiara Razzetta}
\address[C.~Razzetta]{Dipartimento di Matematica, University of Genova, via Dodecaneso 35, IT-16146 Genova, Italy.}
\email{razzetta@dima.unige.it (Corresponding Author)}

\author[V.~Candiani]{Valentina Candiani}
\address[V.~Candiani]{Dipartimento di Matematica, University of Genova, via Dodecaneso 35, IT-16146 Genova, Italy.}
\email{candiani@dima.unige.it}

\author[M.~Crocco]{Marco Crocco}
\address[M.~Crocco]{Esaote S.p.A., Via E. Melen 77, IT-16152 Genova, Italy}
\email{marco.crocco@esaote.com}

\author[F.~Benvenuto]{Federico Benvenuto}
\address[F.~Benvenuto]{Dipartimento di Matematica, University of Genova, via Dodecaneso 35, IT-16146 Genova, Italy.}
\email{benvenuto@dima.unige.it}

\begin{abstract}
 
Medical ultrasound imaging is the most widespread real-time non-invasive imaging system and its formulation comprises signal transmission, signal reception, and image formation. 
Ultrasound signal transmission modelling has been formalized over the years through different approaches by exploiting the physics of the associated wave problem.
This work proposes a novel computational framework for modelling the ultrasound signal transmission step in the time-frequency domain for a linear-array probe.
More specifically, from the impulse response theory defined in the time domain, we derived a parametric model in the corresponding frequency domain, with appropriate approximations for the narrowband case.  
To validate the model, we implemented a numerical simulator and tested it with synthetic data. 
Numerical experiments demonstrate that the proposed model is computationally feasible, efficient, and compatible with realistic measurements and existing state-of-the-art simulators. 
The formulated model can be employed for analyzing how the involved parameters affect the generated beam pattern, and ultimately for optimizing measurement settings in an automatic and systematic way.
\end{abstract}
	
\maketitle
	
\KEY
Medical Ultrasound Imaging, 
Beamforming,
Beam Pattern,
Manifold Domain,
Simulation Software.
\EndKEY
	
\MSC
92C55, 
94A08, 
68U10, 
94A12. 
\EndMSC

\section{Introduction}
Ultrasound imaging is a prominent diagnosis technique, currently used in everyday clinical practice thanks to its moderate cost, radiation-free nature, and the availability of real-time results. Ultrasound imaging, however, suffers from several artifacts such as clutter and poor resolution, originating from its wide Point Spread Function (PSF). 
Most common medical applications include obstetrics, cardiology and surgical guidance \cite{Fenster2015}, where ultrasound imaging is employed to retrieve not only morphological information but it is as well suited for functional imaging.

Commonly used commercial ultrasound scanners gather 2D images, which are created from multiple scans. A standard ultrasound probe consists of multiple transducer elements, used to transmit short acoustic pulses into the considered medium, focused at different directions (transmission step). Tissue inhomogeneities reflect back the acoustic signals, which are received by the array elements. The received signals are then sampled and digitally beamformed to yield a line in the image (receive step). This process is repeated for consecutive directions to create the complete image frame (image formation).

Over the last decades, several mathematical models have been proposed to deal with medical ultrasound beamforming process, which comprises 1) signal transmission, 2) signal reception, and 3) image formation, tackling the modelling problem by means of different approaches.

Papers on time domain modelling of ultrasound transmission and receiving steps \cite{Jensen1991,Jensen2002, Szabo2004} use spatial impulse response techniques, whereas papers working in Fourier domains \cite{Crombie1997, Li1999, Eldar2014, Guo2022, Yu2022} deal with the problem directly in the frequency domain, by considering the spectrum of the signal instead of the waveform in time, with the so-called Fresnel approximation \cite{Greenleaf1994} or from a compressed sensing viewpoint. Both time and frequency domain models then usually employ the Delay-And-Sum (DAS) algorithm \cite{Friis1937} for the actual image reconstruction step. 
Concerning the image formation process as a whole, several papers formulate the ultrasound image reconstruction problem as an inverse problem \cite{Basarab2016,Ozkan2018,Basarab2022}, where they perform beamforming through a regularized inverse problem based on linear models relating the reflected echoes to the to-be-recovered signal. 
Another line of research has been devoted to mitigating the adverse effect of non-ideal PSF using a deconvolution approach \cite{Hundt1980, Michailovich2007, Chen2016}.  
In \cite{Besson2019} the authors revisit standard image formation techniques through the mathematics of an elliptic Radon transform. Finally, deep learning techniques are also strongly emerging in the last few years \cite{Luchies2018,Hyun2019}. 

Existing numerical simulators currently used in medical applications exploit some of the aforementioned models. Matlab simulators comprise  \cite{Jensen1996,kwave2010,Jensen2014, GARCIA2022}, among which the updated version of FIELD II \cite{Jensen2014} has been considered over the years as an absolute benchmark for all medical ultrasound imaging numerical simulations, being the most extensive and comprehensive software. Among Python simulators, \cite{Villalba2019,Cueto2022} can be cited among the most general models, along with more application-specific Github packages and the Python wrappers of the abovementioned models.

This work deals with the formulation, implementation and validation of a numerical model and the corresponding simulator parUST (parallel parametric UltraSound Transmission software) for the signal transmission step only \cite{parust}.
The new formulation addresses the computational problem in a hybrid framework by splitting the computations into two parts: in the first one the transmission step is performed in the time domain while, in the following second part, it is computed in the frequency domain.
This strategy allows us to get results comparable with the state of the art and at the same time greater flexibility in providing beam patterns from different sets of parameters.
Indeed, numerical evidence shows that the beam patterns derived from our novel simulator are essentially comparable qualitatively and quantitatively to the outputs of the state-of-the-art simulator FIELD II \cite{Jensen2014}. 

Our framework has several attractive qualities. Its frequency domain approach avoids time oversampling, necessary in time domain for accurate signal delay computation, thus achieving comparable accuracy while working close to the Nyquist limit, with a consequent saving of memory and computation time, albeit only considering the signal transmission computation.
It also specializes the narrowband case with a single frequency computation, hence allowing an even higher computational saving with respect to FIELD II-like time-domain approaches, in which large time windows needed for narrowband signals are particularly cumbersome. 
Finally, this hybrid approach computes the transmitted signal in the form of a parametric imaging model, thus explicitly extracting and better highlighting its dependency on the model varying parameters, such as the transmit delays, the emitted pulse frequencies, and the considered apertures.  
This allows one to govern the free model parameters, paving the way for a number of different research directions. 
This work was, in fact, motivated by the Transmit Beam Pattern (TBP) optimization \cite{Razzetta2023}, and by the need for the PSF optimization, with the leading aim of ultimately enhancing both the image quality and the quantitative measurements for functional imaging, by generating optimal beam patterns and narrowing down the system PSF. 
Moreover, by considering the model as a parametric imaging problem, one also achieves a deeper understanding of the parameters' behaviour and their domain shape.
For instance, when only $N$ transmit delays are considered as variables and the parameters' domain is an $N$-dimensional torus, we are able to guide the choice of ad-hoc optimization algorithms that exploit the particular topology of the parameter space, such as algorithms designed for manifold optimization and stochastic approaches.  

This paper is organized as follows. Section 2 reviews the model formulation in time domain and proposes its frequency domain counterpart in both the wideband and the narrowband case. 
The model implementation is presented in Section \ref{sec:software} along with the approximations introduced to speed up the computations.
Section 4 exhibits the simulation setting and the numerical experiments, and conclusions are drawn in Section 5.

\section{The transmission model formulation}\label{sec:model}
Medical ultrasound systems can be described through the physics of emitting and receiving acoustic waves, as well as their interaction with the surrounding material.
As the beamforming process theory has been extensively studied, in this section we will focus on describing a hybrid model in time-frequency domain exclusively for the transmission part.
The aim is to analyze how the transmission energy pattern varies according to all the involved parameters. To begin with, the more general wideband case will be analyzed, followed by the more specific narrowband case and its related approximations.

\subsection{Wideband modelling}\label{ssec:wide}
\begin{figure}
    \centering
    \includegraphics[width=0.8\textwidth]{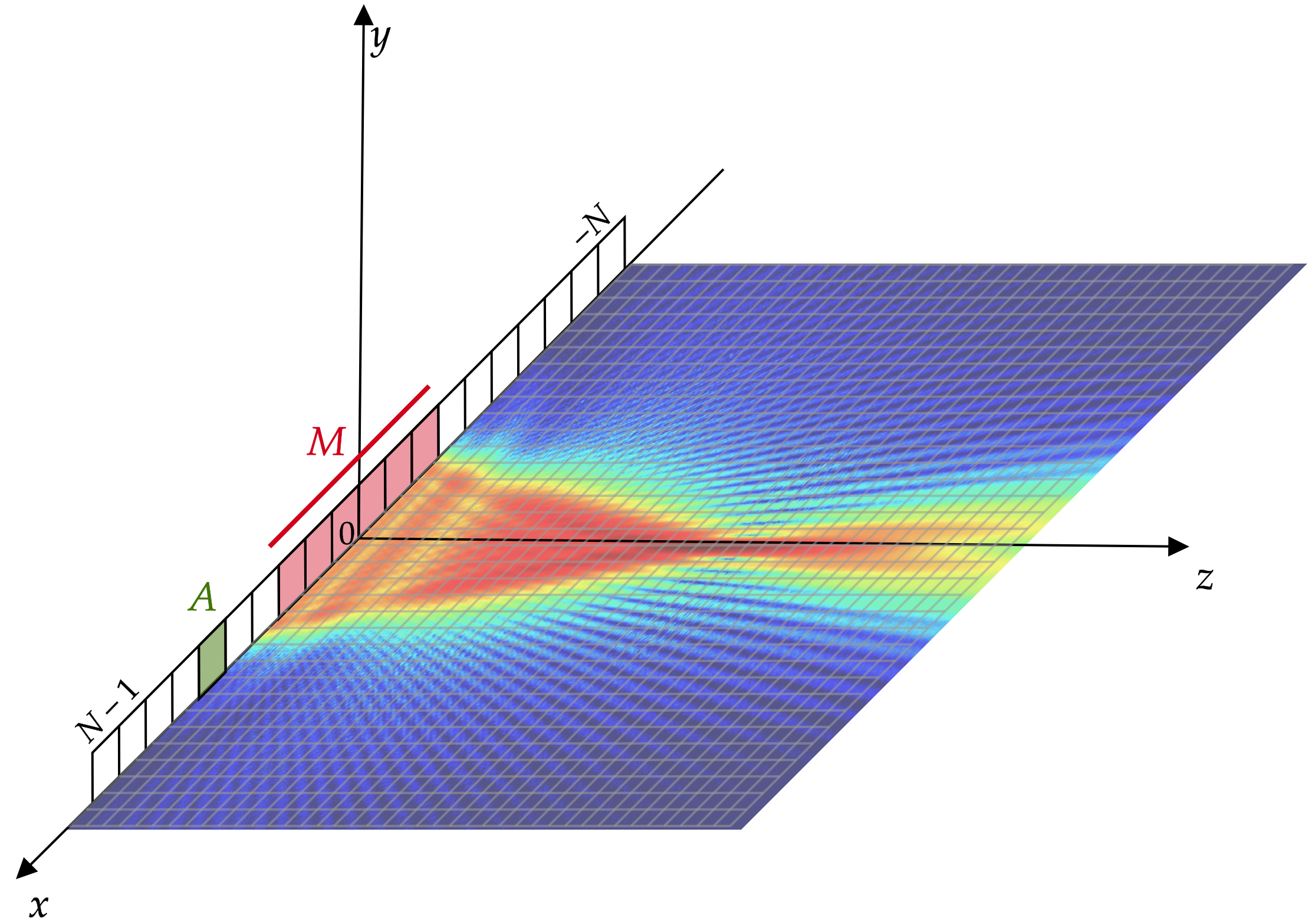}
    \caption{Depiction of a beam pattern resulting from a linear-array probe with $2N$ elements of area $A$, of which only $M$ are activated, symmetrically with respect to the probe central axis.}
    \label{fig:schema}
\end{figure}

Each electric linear system can be identified by its impulse response function, which describes the way a system reacts to the emission of a Dirac pulse.
Let us fix a coordinate system such that the probe is placed on the $xy$ plane, while the field space is given by $\Omega := \left\{ \Vec{r} := (x,y,z) \in \mathbb{R}^3 \vert z > 0 \right\}$, as depicted in Figure \ref{fig:schema}. The impulse response at each field point depends on the distance between the probe and the point itself \cite{Harris1981}.
We consider a probe of surface $S$ emitting a pulse with speed $v(\Vec{s},t)$, outward-pointing and orthogonal to the probe surface, dependent on the position on the probe, and a target point $\Vec{r}$ in the field. At each temporal instant $t$ the impulse response is given by:
\begin{equation}
    h_S(\Vec{r}, t) =\int_S \frac{v\left(\Vec{s}, t - \frac{\lVert \Vec{r}-\Vec{s}\rVert}{c}\right)}{2\pi\lVert \Vec{r}-\Vec{s}\rVert} \, d\Vec{s},
    \label{eq:h}
\end{equation}
where $c$ is the medium speed of sound in tissues.

Medical ultrasound arrays are composed of several piezoelectric elements.
These elements are geometrically identical and it is reasonable to assume that they emit the same waveform, independently of the position on the probe \cite{Jensen1996}.
Thus, we consider a probe of $2N$ elements of width $\Delta x$ and surface $A$.
The impulse response function for a fixed element 
depends only on the mutual distance between the element and the point.
A transmitted Beam Pattern (BP) corresponds to the total energy flux in the field over time.
In our setting the transmission step is performed by sequentially activating several elements of the probe emitting the same waveform.
For each $n \in \mathcal{N} = \{-N, \dots, N-1\}$, let an element be a rectangle with area $A_n$, centered at position $\Vec{c_n} := (n+1/2, 0, 0) \Delta x$ and fixed $\Delta y$.
We then consider only a subset of active elements that we indicate with $M \subseteq \mathcal{N}$ and $m \in M$.
If we denote by $I(t)$ the emitted waveform and by $a(\Vec{r},t)$ the thermoviscous wave attenuation in a homogeneous medium \cite{Jensen1999,Szabo2004}, we get:
\begin{equation}
    H_m(\Vec{r},\cdot) := a(\Vec{r},\cdot) \ast h_{A_m}(\Vec{r},\cdot) ~ .
    \label{eq:HI}
\end{equation}

While crossing the field, the waves will sum up coherently and incoherently generating different patterns.
It holds that the signal crossing a point at a fixed time instant is:
\begin{equation}
    S(\Vec{r}, \cdot) = I \ast \sum_{m\in M} H_m(\Vec{r}, \cdot) ~ .
    \label{eq:Signal}
\end{equation}
Accordingly, the energy flux at each point can be computed as the integral over time of the square module of such signal:
\begin{equation}
    E(\Vec{r}) = \int_{\mathbb{R}} \vert S(\Vec{r}, t)\vert ^2 \, dt.
    \label{eq:E}
\end{equation}
The simultaneous emission of multiple identical pulses generates a plane wave that propagates through the field.
Even though plane wave transmission and receiving modalities have been developed \cite{Tanter2014}, there is still a massive use of the so-called focused transmissions.
The focusing process consists in emitting the waves at different times, so that the resulting energy will be higher at certain field points due to coherent sums \cite{Cobbold2007}.
These emission time differences between elements can be described as a vector $D=(D_1,\ldots,D_M)$ of shifting in time of the signals, so that the delayed signal at point $\Vec{r}$ is:
\begin{equation}
    S_D(\Vec{r}, \cdot) = I \ast \sum_{m\in M} H_m(\Vec{r}, \cdot) \ast \delta(t-D_m).
    \label{eq:Sdel}
\end{equation}
Consequently, the total energy flux at $\Vec{r}$ is obtained as:
\begin{equation}
     E(\Vec{r}) = \int_{\mathbb{R}} \vert S_D  (\Vec{r}, t)\vert ^2 \, dt.
    \label{eq:Edel}
\end{equation}
We can observe that the relevant time interval is determined by the time at which the pulse actually crosses the point, i.e., the emitted wave length.
To this end, we consider the Fourier transform of Eq. \eqref{eq:Sdel}:
\begin{equation}
    \hat{S}_D(\Vec{r}, f) = \hat{I}(f)\sum_{m\in M} \hat{H}_m(\Vec{r}, f)e^{- 2 \pi i D_m f},
    \label{eq:Sfreq}
\end{equation}
where $\hat{I}(f)$, $\hat{H}_m(\Vec{r}, f)$, and $e^{- 2 \pi i D_m f}$ are the Fourier transforms of the emitted waveform, the attenuated impulse response function, and of the time delay, respectively.

By Plancherel's theorem, it holds:
\begin{equation}
    E(\Vec{r}) = \int_{\mathbb{R}} \vert S_D  (\Vec{r}, t)\vert ^2 \, d t = \int_{\mathbb{R}} \vert \hat{S}_D  (\Vec{r}, f)\vert ^2 \, df.
    \label{eq:Efreq}
\end{equation}
 The energy flux depends on the waveform, on the number of emitted pulses, and on the delays applied for the focusing step.
This formulation is valid for any choice of emitted waveform and transmission pattern adopted by medical ultrasound apparati, on imaging and quantitative modalities. 
For the sake of simplicity we neglected in the model the transmit apodization, i.e. the scalar weighting of the transmit waveform as a function of the elements, and the possibility to transmit different waveforms at each element. These features can only be found on top-level commercial scanners or research scanners. 
If necessary, Eq. \eqref{eq:Sfreq} can be easily updated with the insertion of a scalar multiplicative factor, function of $m$ in the sum for apodization, or by making $\hat{I}(f)$ dependent on $m$ and including it in the sum for element-dependent waveforms.

To simplify the computations, it is possible to restrict the frequency domain to the band of the probe, since frequency components of the transmit waveform outside the band are filtered out.  
In the next paragraph, we specialize the analysis to the narrowband pulse case, thus exploiting some reasonable approximations.

\subsection{Narrowband case approximation}\label{ssec:narrow}
\begin{figure}[!h]
    \centering
    \includegraphics[width=0.8\textwidth]{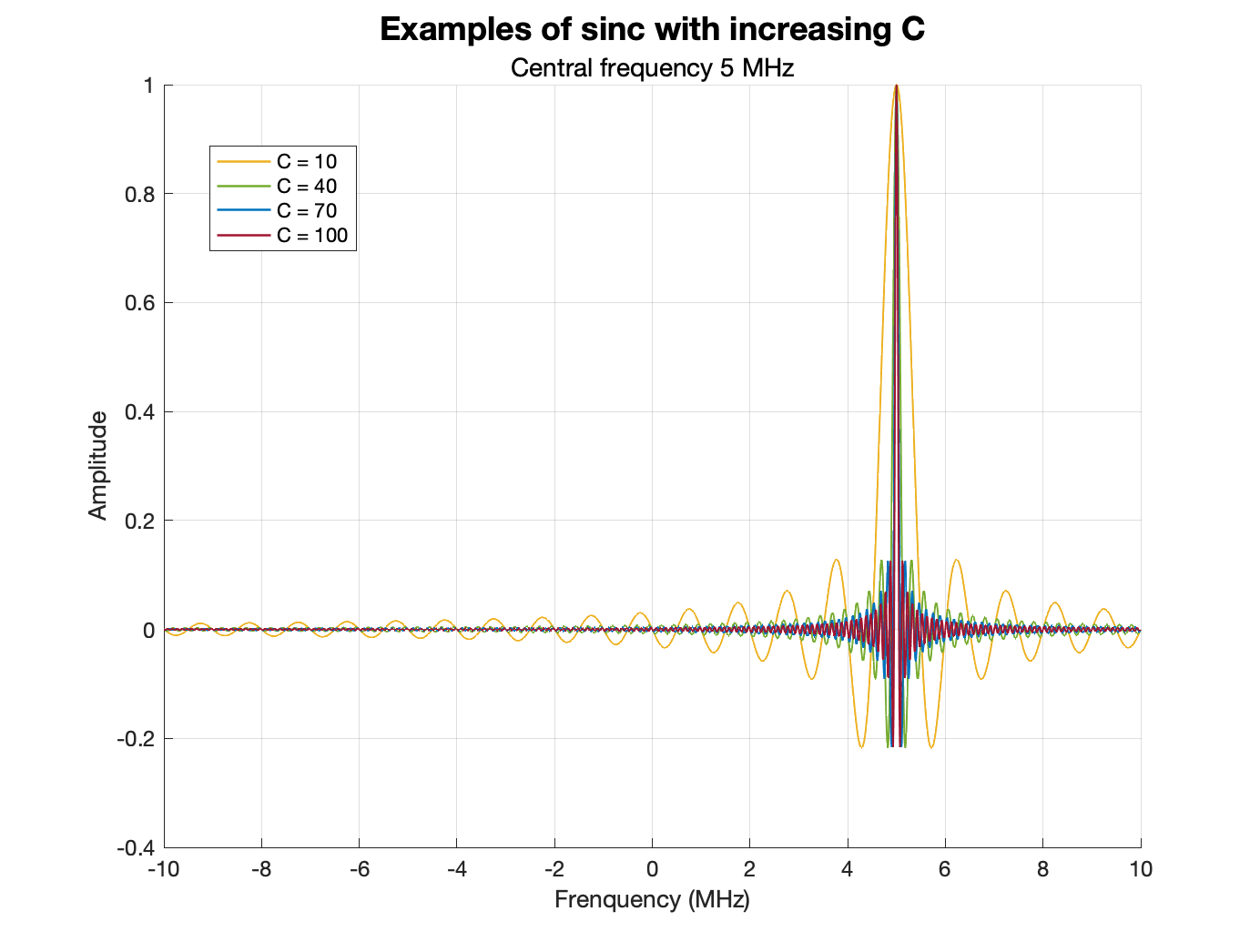}
    \caption{Behaviour of the sinc function with increasing number of cycles $C$.}
    \label{fig:sincPlot}
\end{figure}

While the large majority of ultrasound modalities deals with wideband beam patterns, two specific quantitative applications, Continuous Wave (CW) Doppler \cite{Brody1974} and Acoustic Radiation Force Impulse (ARFI) elastography \cite{ARFI3} employ narrowband beam patterns.
Assuming the emitted waveform to have a narrow band is equivalent to considering that during the emission time many cycles of the pulse are present.
The high number of cycles makes it physically reasonable to approximate the waveform as an infinite pulse, with the following rationale.
One can distinguish between different types of transmissions, by analyzing the fractionary band of the emitted signal.
\[ b_f : = \frac{f_{max}- f_{min}}{f_0},\]
with $f_0$ being the central frequency of the pulse.
In some quantitative applications of medical ultrasound, such as ARFI Elastography, the transmission consists of hundreds of cycles, i.e., a very small fractional band.
In this setting, as the number of pulse cycles increases, the resulting beam pattern can be approximated fairly accurately with a beam pattern related to a single frequency time infinite signal.
As a matter of fact, the tighter the band the lesser frequencies are significant for the computation of the energy flux (cf. Eq. \eqref{eq:Efreq}).
As a limit case, the overall signal described in Eq. \eqref{eq:Sfreq} differs from zero just at the central frequency.
To formalize this phenomenon, we assume the waveform to be uniquely characterized by the pulse central frequency $f_0$.
Hence, an emitted pulse  $\Tilde{I}$ can be described as an infinite complex exponential $I$ multiplied by a rectangular window $\Pi_a$:
\begin{equation}
\Pi_a\left(t\right) \coloneqq \begin{cases}
    0 \quad \text{if} \quad \vert t \vert \geq \frac{a}{2} \\
    1 \quad \text{if} \quad \vert t \vert < \frac{a}{2}
\end{cases}
\end{equation}
whose length corresponds to the emission time interval $\tau$:
\begin{equation}
    \Tilde{I}(t) = I(t) \, \Pi_\tau\left(t\right) = e^{-2\pi i f_0 t} \, \Pi_\tau\left(t\right).
    \label{eq:pulse}
\end{equation}
Switching to its Fourier transform, we obtain the cardinal sinus centered at $f_0$:
\begin{eqnarray}    
    \hat{\Tilde{I}}(f) &=& \left[ \mathscr{F}(I(t)) \ast \mathscr{F}\left(\Pi_\tau\left(t\right)\right) 
    \right](f) \\ \nonumber
    &=&  \left[ \delta(\cdot-f_0) \ast \operatorname{sinc}\left(\tau \cdot\right) 
    \right](f)= \operatorname{sinc}\left(\tau (f-f_0)\right).
    \label{eq:Ifreq}
\end{eqnarray}
It should also be observed that the considered time interval contains a finite number $C$ of pulse cycles of frequency $f_0$:
\[\tau = \frac{C}{f_0},\]
so that for $C \gg 0$, it holds that $\operatorname{sinc}\left(\tau (f-f_0)\right) \approx \delta(f-f_0)$, since $\tau \gg 0$. Such behavior can be appreciated in Figure \ref{fig:sincPlot}.
As a consequence, for a high number of cycles, i.e. narrowband signals, it holds:
\[\hat{\Tilde{I}}(f) \approx \delta(f - f_0).\]
Then, the delayed narrowband signal of Eq. \eqref{eq:Sfreq} ends up being:
\begin{equation}
    \hat{S}^{f_0}_N(\Vec{r}) : = \delta(f-f_0)\sum_{m\in M} \hat{H}_m(\Vec{r}, f)e^{- 2 \pi i D_m f}.
\end{equation}
It directly follows that in this setting the energy flux cannot be defined as the integral over the frequencies, the form being divergent.
We can consider instead the power of the signal:
\begin{equation}
    P(\Vec{r}) = \left|\sum_{m\in M} \hat{H}^{f_0}_m(\Vec{r})e^{- 2 \pi i D_m f_0}\right|^2,
    \label{eq:Power}
\end{equation}
where $\hat{H}^{f_0}_m(\Vec{r}) : = \hat{H}_m(\Vec{r}, f_0)$ is the evaluation of the attenuated impulse response Fourier Transform at frequency $f_0$.
It occurs that the beam pattern, seen as a collection of the power values, depends on the central frequency $f_0$, the number of active elements $M$ and the delay values $D_m$, for $m=1,\dots, M$, which from here on we will consider the model free parameters.
\subsubsection{Symmetry with respect to the center}
We will focus on the (realistic) case in which the transmission is performed with an equal number of active elements with respect to the central transmission axis (i.e., the geometrical center of the probe, invariant under translations on the $x$ axis) and without beam steering.
In this setting, symmetrical elements will have the same impulse response and Eq. \eqref{eq:Power} can then be rewritten as:
\begin{align}
    P(\Vec{r}) &= \left(\sum_{m\in M} \hat{H}^{f_0}_m(\Vec{r})e^{- 2 \pi i D_m f_0}\right) \overline{\left(\sum_{n\in M} \hat{H}^{f_0}_n(\Vec{r})e^{- 2 \pi i D_n f_0}\right)} \nonumber\\
    &= \sum_{m,n \in M} \hat{H}^{f_0}_m(\Vec{r})\overline{\hat{H}^{f_0}_n(\Vec{r})}e^{- 2 \pi i (D_m-D_n) f_0}\nonumber\\
    &=2 \sum_{\substack{m\ge n \\ m,n \in \{0, \dots, \Bar{M}\}}} \hat{H}^{f_0}_{m,n}(\Vec{r}) \cos(2\pi f_0 (D_m -D_n)),
\end{align}
where $\Bar{M}$ is the maximum positive integer index of the symmetric aperture, and we define
\[ \hat{H}^{f_0}_{m,n}(\Vec{r}) : = \hat{H}^{f_0}_m(\Vec{r})\overline{\hat{H}^{f_0}_n(\Vec{r})} = \hat{H}^{f_0}_n(\Vec{r})\overline{\hat{H}^{f_0}_m(\Vec{r})}.\]


\subsection{Interpretation of the parameter domain as a manifold}
In the previously described symmetrical setting and under the assumption of applying the same delay to symmetrical elements (reasonable for all the settings in which steering is not needed), we can consider the power of the signal as a parametric model $ \mathscr{P}$ depending on the parameters $f_0, \, D_0, \dots, D_{\Bar{M}}$:
\begin{equation}
    P(\Vec{r}) = \left|\sum_{m \in \{0, \cdots, \Bar{M}\}} \left(\hat{H}^{f_0}_{m} + \hat{H}^{f_0}_{-m-1} \right) e^{-2\pi i D_m f_0} \right|^2 =: \mathscr{P}(f_0, D_0, \dots, D_{\Bar{M}}),
    \label{eq:Psim}
\end{equation}
where the exponential term represents a phase shift for the corresponding sum of impulse responses, and $\Bar{M}$ replaces $M$ in the symmetrical setting. 
\begin{figure}[ht]
    \centering
    \includegraphics[width=0.9\textwidth]{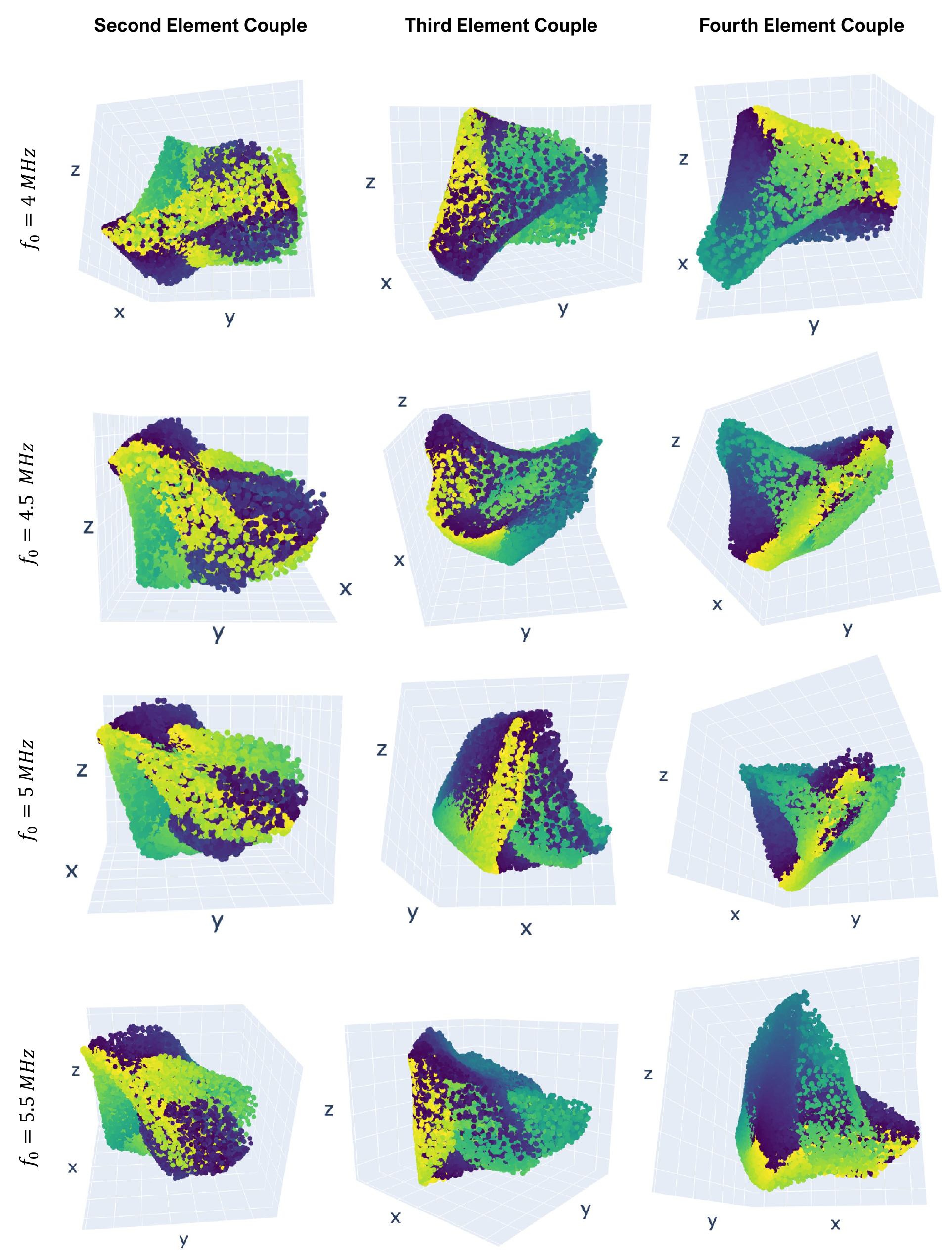}
    \caption{Examples of the parameters' domain for a linear-array probe. From the top, each row corresponds to a different fixed frequency $f_0 = 4e6, 4.5e6, 5e6, 5.5e6 \, MHz$, where $M = 8$ ($\Bar{M}=4$) elements are active, with the two central ones not delayed.}
    \label{fig:tori}
\end{figure}
It is easy to see that, for a fixed frequency $f_0$, referring to $D_m$ as the variable, the domain of the exponential function is a circumference.
Thus, for a fixed frequency $f_0$, for $\bar{M}$ active elements, and under the assumption of independent delays, the domain of the parametric model form will contain an $\Bar{M}$-dimensional torus:
\begin{equation}
    \mathbb{T} = \underbrace{S^1\times \cdots \times S^1}_{\bar{M} \text{times}}.
    \label{eq:torus}
\end{equation}
Figure \ref{fig:tori} displays different parameter domains for the same linear-array probe setting, with varying central frequency of the narrowband pulse along the rows and  with eight active elements, six of which are delayed symmetrically with respect to the center. 
In this case the free parameter space is a 4-dimensional torus.

For each frequency we have computed $30\,000$ BPs, with randomly symmetric delays.
By performing a Principal Component Analysis we have projected the sets of BPs from the pixel space onto a $3$D space to visualize the $4^{\text{th}}$ dimension manifold.
Each column corresponds to a delay in that the colour represents the phase shift: blue represents the minimum angle $0$ and yellow represents the maximum angle $2\pi$.
Moreover, the yellow and blue dots are always close as in euclidean spaces the two angles are different points in the space, while in our domain they correspond to the same point.
This interpretation is the starting point for the TBP optimization explained in \cite{Razzetta2023}, where we propose to optimize the narrowband BP by considering each single transmit delay as an independent variable and keeping the transmit intensity fixed.
This understanding will also lay the foundations for future works on automatic PSF optimization. 

\section{Software implementation}\label{sec:software}
To validate the presented transmission model, we have implemented the high-performance numerical simulator parUST running on Python and freely available on GitHub \cite{parust}. 
The simulator is currently tailored for linear-array probes, but its implementation is easily extendable to convex ones with a trivial change of coordinates from Cartesian to polar.
We tested the simulator against FIELD II \cite{Jensen2014} and compared the results in terms of generated beam patterns and computational efficiency in terms of running time. 
It should be noted, however, that we do not hereby claim the superiority of our algorithm with respect to the more complete and extensive FIELD II, but we supply an open and faster alternative exclusively for the transmission step. 

Since the computational cost depends significantly on the computation of the impulse response of the probe elements, we have structured the simulator in such a way that, for a defined probe model and field region of interest, these computations are performed only once.
Hereafter, we will call \textit{map} the ordered collection of all the impulse responses for a fixed couple of field grid and probe element.
We proceed with describing our simulation tool by dividing two steps: 1) the map generation and 2) the BP generation.
At each step, we need to distinguish between the implementation for wideband and narrowband cases since the mono-frequency approximation in the model determines significant differences in the size of arrays and memory management. 

Concerning the map generation, we provide a common implementation for the impulse response function calculation and for the probe description used by both wide and narrow settings.
The computations are structured to be performed just once and offline for each setting by parallelizing as much as possible and thus exploiting the available cores without involving GPU computation \cite{Bruyneel2013}. 
We need to estimate the value of the impulse response in time for each couple (element of the probe - point of the field).
Referring to the definition in Eq. \eqref{eq:HI}, for each point one needs to compute the integral over the element surface of a spherical wave propagating from the point itself, as it is done in FIELD II \cite{Jensen1996, Jensen2014} with the approach described in \cite{Jensen1999H}. We use, however, a different technique.
The following implementation choices lay the foundations for a fast, flexible, and scalable open-source Python package.

\subsection{Field discretization and probe description}
We assume the field to be on the $xz$ plane at $y=0$ and we discretize it by sampling the space such that the $x$ coordinates are taken as sub-multiples of the element pitch, the sum of the element width, and the space between \textit{kerf} two consecutive elements (cf. Figure \ref{fig:hcalc}).
Being $\Delta x$ the element pitch and $\Delta z$ the discretization step along $z$, the grid points are defined as follows:
$$\{(x_u, z_v) : u = -LN, \cdots, L(N-1), v = 0, \cdots, N_Z, L \in \mathbf{Z}_+\},$$ 
$$\text{where} \,\ x_u = \frac{u\Delta x}{L} \,\ \text{and} \,\  z_u = \frac{v\Delta z}{N_Z}.$$

The probe is composed of several identical elements, whose surface is discretized 
by adopting a grid of $Nx$ points along the $x$ axis and $Ny$ points along the $y$ axis, indexed with an unrolled index $j$ such that:
$$j = l_1 + N_x l_2, \,\ \text{where} \quad l_1 = -N_x/2, \cdots, N_x/2 \quad \text{and} \quad  l_2 = -N_y/2, \cdots, N_y/2.$$
The probe surface is usually covered with an acoustic lens that focuses the signals over the image plane $xz$.
In detail, each probe lens has its own focus point, called 
\textit{geometrical focus} $\Vec{g}$ $$\Vec{g} = (0,0,z_g).$$
We describe this characteristic by associating to each point of the element discretization a time delay that encodes the differences in time of flight between the geometrical focus and the point on the element:
$$\tau_j^g = \frac{\sqrt{y_j^2+z_g^2}}{c},$$
where the $x$ coordinate is not taken into account since the effect of the lens is identical for each choice of azimuth and it modifies just the pulse in the elevation direction.
\subsection{Computation of \textit{h}}
\begin{figure}
    \centering
\includegraphics[width=0.65\textwidth]{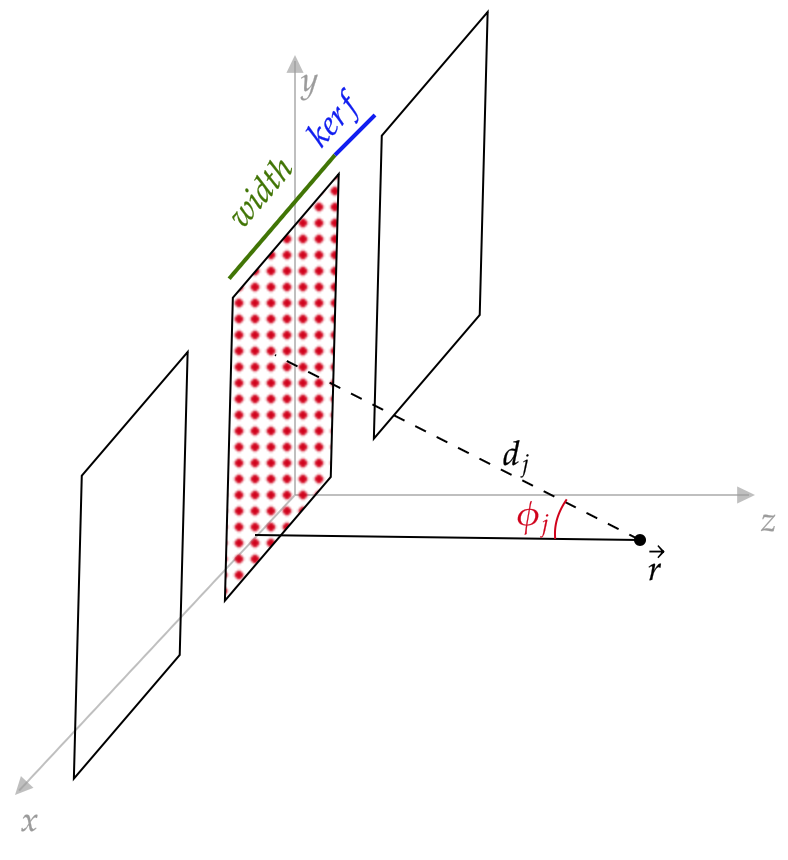}
    \caption{Depiction of the geometry adopted for the calculation of a spatial impulse response for a fixed field point $\Vec{r}$ and a probe element. The element discretization grid is represented by the red dots, each of which is associated with the distance $d_j$ and the angle $\phi_j$.}
    \label{fig:hcalc}
\end{figure}

We can consider each point of the element discretization as a small source of a spherical wave (Figure \ref{fig:hcalc}) and approximate the impulse response as the sum of contributes that arrive to a field point at a specific time instant.
The sum of all spherical waves must correspond to the wave emitted by the entire element.
This can be achieved by delaying and weighting the single wave taking into account the distance between the source and the field point $d_j$.
Thus, if we consider the soft baffle apodization defined in \cite{Jensen1999} and introduce $\phi_j$ as the angle between the minimum trajectory from the point to the probe plane $z = 0$ and the trajectory between the point and the probe point $j$, the single wave is:
$$w_j(t) = \frac{\delta(t - \tau(d_j))}{d_j} \cos(\phi_j) ~ .$$

In order to implement the computation on a discrete time scale we need to limit the frequency band of the Dirac function before time sampling so as to avoid aliasing. To this end, the Dirac function is convolved with a cardinal sinus, whose Fourier transform is a rectangle centered at zero and of width equal to $f_{s}$, with $f_s$ being the sampling frequency constrained to be higher than twice the maximum frequency of the probe band:
$$w_j(t) \approx \frac{\operatorname{sinc}(f_s(t-\tau(d_j)))}{d_j}\cos(\phi_j).$$
For a fixed field point the computation of $h$ is performed as described in Algorithm \ref{alg:ALG1}.
\begin{algorithm}
\caption{Computation of $h$ for $\Vec{r} $}
\begin{algorithmic}
\Require
\State $d$ all the distances between the element discretization points and $\Vec{r}$
\State $\tau = d/c-\tau^g$ all the times of flight shifted by the correct geometrical delay
\For{$j$ in $(0, J)$} \do \\
\State $h(t) = h(t)  + w_j(t)$
\EndFor
\end{algorithmic}
\label{alg:ALG1}
\end{algorithm}
The code for Algorithm \ref{alg:ALG1} is in c++ language and it is parallelized to use all available CPUs for the calculation over all field points.
Finally, the collection of impulse responses is processed by the Fast Fourier Transform (FFT) along the time axis.
From this point on, all calculations will be performed in time frequency domain to reduce the computational cost.

\subsection{Wideband maps}
In the wideband case we need a map for each active element of the probe, with dimensions (field grid) $\times$ (number of computed frequencies).
This high dimensional matrix may be too big to be stored or saved for future reuse.
For this reason we proceed with the following steps:
\begin{enumerate}
    \item compute one impulse response map containing all the impulse responses of the field points with respect to a central probe element,
    \item save the map,
    \item exploit the translation invariance along the probe: identical elements will have the same impulse response, only spatially translated by a pitch.
\end{enumerate}
Since we are dealing with real-valued signals we can limit the map computation to the first half of the spectrum, corresponding to positive frequencies, hence halving the dimension of the maps.
\subsection{Narrowband maps}
As for the narrowband case, we adopt the single frequency approximation as described in Section \ref{ssec:wide}, so that the memory usage and the computational cost are significantly reduced.
In this case, after computing the map for the central frequency, we proceed by exploiting the central axis symmetry and by saving the maps summed in pairs.
The choice of computing the entire maps and save them is coherent with the following remark: the computation is time consuming even if parallelized, but it depends solely on the probe characteristics and on the grid dimensions.
Once the computation is handled, it is possible to perform different experiments with the same basic settings and only tuning the number of active elements, the delays and the emitted pulse.
\subsection{Beam pattern generation}
After computing the maps, the beam pattern calculation is straightforward and considerably fast (multiplication of pre-computed matrices and arrays), with the wideband case being slightly more computationally expensive due to the number of considered frequency.
To reduce the computational cost we allowed taking into account the probe band in the map calculation step.
The advantage is to consider only the few significant frequencies that can be registered by the probe, namely discarding the ones that differ less than $40 \, dB$ from the highest one.

\section{Numerical results}\label{sec:results}
This section presents our numerical tests.
More specifically, we compare qualitatively and quantitatively the same experiments performed with FIELD II \cite{Jensen2014} (in Matlab) and our Python simulation tool, with the necessary adjustments in switching from one programming language to the other.
For both wideband and narrowband cases we consider a linear-array probe composed of $192$ elements of height $5 \, mm$ and pitch $0.245 \, mm$, with geometrical focus at $25 \,mm$.
The field is sampled with step $dx = \frac{pitch}{4}$ along the $x$ axis and $dz = 0.2 \, mm$ along the $z$ axis by taking into account a minimum depth of $2 \, mm$ and a maximum depth of $42 \, mm$.
The adopted time step of $dt = 10^{-8}$ corresponds to a sampling frequency of $100 \, MHz$.
Numerical tests have been performed using two values for the central frequency of the emitted pulse, $f_0 = 3 \, MHz$ and $f_0 = 4.5 \, MHz$, two different numbers of active elements $M = 20$ and $M = 50$, and three depths of focus $F = 10, 25, 35 \, mm$.
\subsection{Quantitative measure}
To quantitatively evaluate the reliability of our simulator parUST in both cases we use a quantitative metric that estimates the minimum distance between the considered pairs of beam patterns.
From now on we will denote by $BP^F$ the collection of values computed by FIELD II and by $BP$ our simulator approximation, each one normalized by its own maximum value, assuming they are computed with the same setting (probe, grid, field...), and thus having the same number of field grid points $N_p$.
As beam patterns are computed up to a normalization constant, we evaluate their distance as follows
\begin{equation}
D(BP, BP^F) \coloneqq \frac{1}{N_p} \min_{\alpha >0} \left\lVert BP - \alpha BP^F\right\rVert^2.
\label{eq:distance_min}
\end{equation} 
Taking the derivative with respect to $\alpha$ to compute the minimum of this convex function, it holds that
\begin{equation}
D(BP, BP^F) = \frac{1}{N_p} \left\lVert BP - \alpha^* BP^F\right\rVert^2,
\label{eq:distance}
\end{equation}
where
\[\alpha^* = \frac{<BP, BP^F>}{<BP^F, BP^F>}.\]
It is worth noticing that such distance depends on the number of pixels chosen for representing the beam patterns, and hence on their spatial resolution: the finer the discretization the larger the error and viceversa. 
On the other hand, the spatial resolution cannot exceed a certain limit due to the physical and geometric properties of the measurements.

\begin{figure}
    \centering
    \includegraphics[width=1.\textwidth]{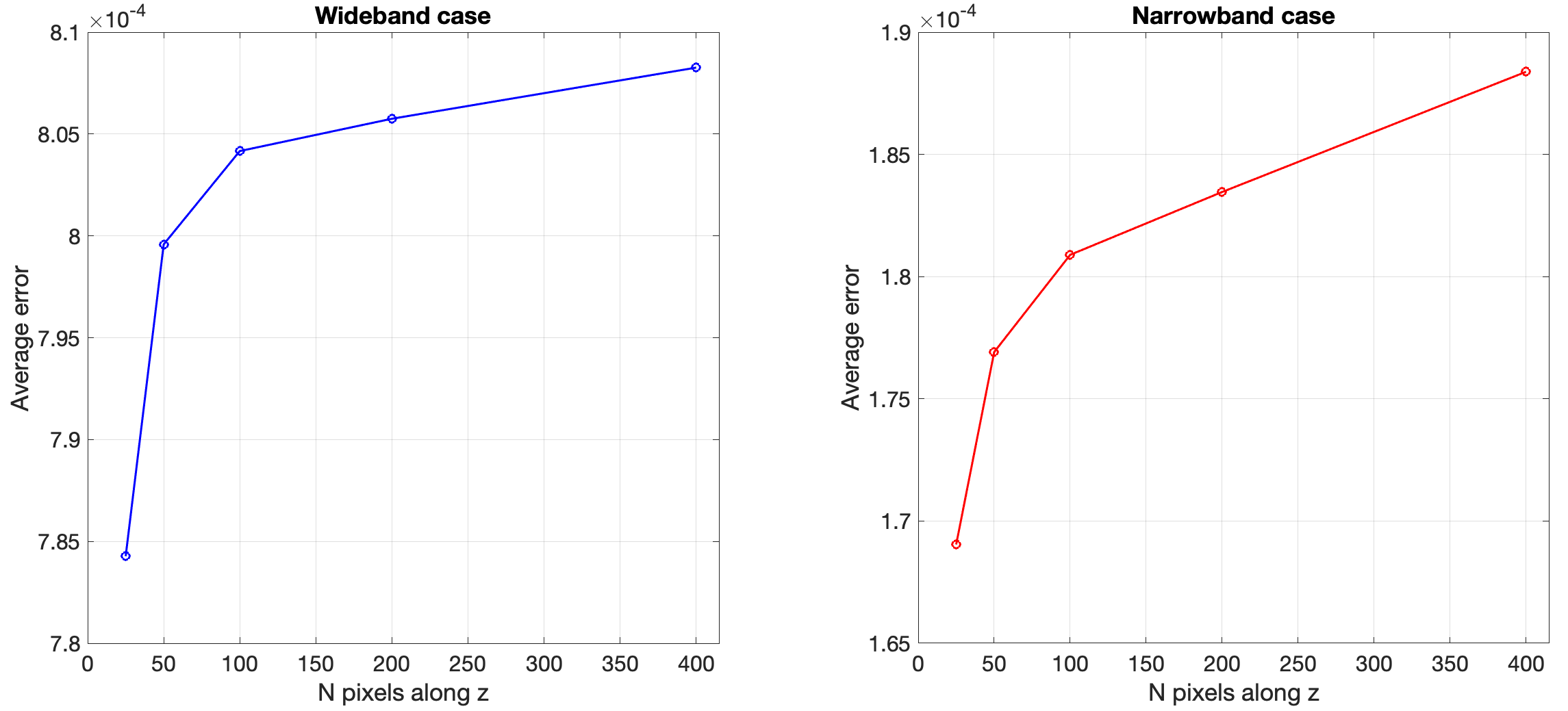}
    \caption{Average BP error $D(BP, BP^F)$ on the $12$ wideband and $12$ narrowband cases, as a function of the number of points along depth. The error is computed for $5$ different resolution sizes along depth, with a fixed $dx = \frac{pitch}{4}$.}
    \label{fig:AverageErr}
\end{figure}

To better justify the beam pattern results, we have computed the error with different spatial resolutions, averaged on all $12$ combinations of values above.
As the average BP error measures the difference between two discretized functions, 
the finer the discretization, the larger the error (cf. Eq. \eqref{eq:distance} and Figure \ref{fig:AverageErr}).
Along the $x$ axis we have fixed the step $dx = \frac{pitch}{4}$ which is the maximum resolution usually achieved in azimuth with ultrasound scanners.
Along depth, we have considered five different resolutions, dividing the $4  \, cm$ field into $25, \, 50, \, 100, \, 200, \, 400$ points, hence reaching a maximum resolution of $0.2 \, mm$.
We can notice that the actual axial resolution of the final image is usually linked to the wavelength of the emitted pulses and the considered minimum discretization step corresponds to a wavelength shorter then the ones related to the transmitted frequencies of $3$ and $4.5 \, MHz$. 
In the following, we present the results obtained when considering $400$ points along depth, i.e., the most accurate measure of BP discrepancy.
\subsection{Test 1: Wideband case, two frequencies, two apertures, three focusing depths}
In the wideband case we have used a sinusoidal pulse of just two cycles, resulting in the computational time for the beam pattern construction being on average $54 \,s$ with our simulator and $80 \,s$ with FIELD II.
From the quantitative results reported in Table \ref{tab:WideQuant}, it can be noted that the average absolute distance $D$ between the couples of simulated signals is $8.08 \cdot 10^{-4}$.

\begin{table}[!ht]
\begin{tabular}{|c| c |c| c |c|}
\cline{2-5} \multicolumn{1}{c|}{} & \multicolumn{2}{c|}{\cellcolor{red!30}Frequency 3.0 MHz}& \multicolumn{2}{c|}{\cellcolor{red!30}Frequency 4.5 MHz} \\ 
\hline
 \multicolumn{1}{|p{1.5cm}|}{\centering \cellcolor{gray!30}F (mm)} & \multicolumn{1}{p{1.5cm}|}{\centering \cellcolor{orange!30}$M = 20$} & \multicolumn{1}{p{1.5cm}|}{\centering\cellcolor{yellow!20}$ M = 50$} & \multicolumn{1}{p{1.5cm}|}{\centering\cellcolor{orange!30}$M = 20$} & \multicolumn{1}{p{1.5cm}|}{\centering\cellcolor{yellow!20}$M = 50$} \\ 
\hline
10 & 4.91e-4 & 9.11e-5 & 2.17e-4& 6.43e-5\\
\hline
25 & 1.05e-3 & 9.20e-5 & 3.57e-4 & 4.25e-4 \\
\hline
35 & 1.21e-3 & 2.83e-3 & 4.51e-4 & 1.58e-3\\
\hline
\end{tabular}
\caption{Quantitative measure $D$ for comparing $12$ experiments in the wideband case. Each row contains the values for a different choice frequency, number of active elements M, and depth of focus F.}
\label{tab:WideQuant}
\end{table}

\begin{figure}
    \centering
    \includegraphics[width=1.\textwidth]{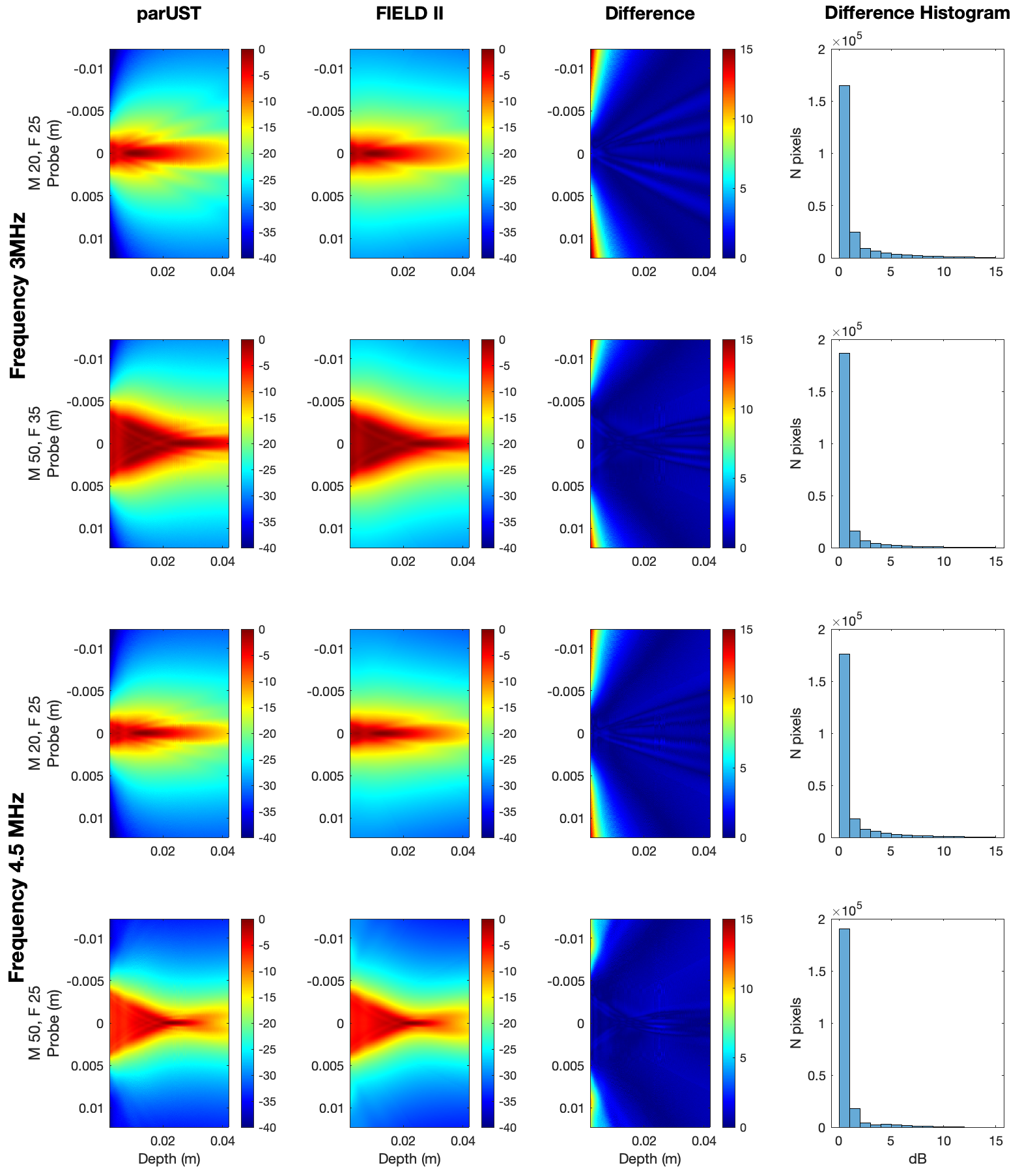}
    \caption{Comparison of generated wideband BPs in different settings. On each column from the left: BPs generated with our simulator parUST, BPs from FIELD II, their absolute difference, and the distribution of pixel value differences. The unity of measure is dB.}
    \label{fig:WidePanel}
\end{figure}

From a qualitative viewpoint, Figure \ref{fig:WidePanel} reports four examples displaying on the first two columns the simulated beam patterns, on the third column their absolute difference and on the fourth column the histogram of the differences. All quantities are expressed in dB. 
The BP shapes are clearly comparable, with almost the same shape in the highest energy region, i.e. the main lobe region.
The major differences (yet negligible in terms of absolute error) are located on the field parts closer to the probe (cf. Figure \ref{fig:WidePanel}, third column). This could be explained with the different approach used to model the element impulse response and in particular the elevation focus, either as a pure delay in our case or through a geometrical shaping of the surface in FIELD case.
Finally, it is immediate to see from the histograms on the fourth column that most pixels differ in intensity less than 5 $dB$, with the highest differences in $dB$ corresponding to the regions of low field intensity.
\subsection{Test 2: Narrowband case, two frequencies, two apertures, three focusing depths}
To simulate the narrowband case in FIELD II we have approximated the single emitted pulse as a sinusoidal wave with a hundred cycles.
In this setting, the need of multiple cycles increases the computational times in FIELD with an average of $2 \, min$, while for our simulator the computation is faster, lasting on average $0.3 \, s$.
From Table \ref{tab:NarrowQuant} we can observe that the errors are lower than in the wideband case with an average of $1.88 \cdot 10^{-4}$.
This is consistent with the differences in computation: the small numerical errors for the impulse response approximation will be propagated through the entire flowof calculations, and the number of operations in the wideband case is larger due to the presence of the entire signal frequency spectrum.
\begin{table}[ht]
\begin{tabular}{|c| c |c| c |c|}
\cline{2-5} \multicolumn{1}{c|}{} & \multicolumn{2}{c|}{\cellcolor{red!30}Frequency 3e6 Hz}& \multicolumn{2}{c|}{\cellcolor{red!30}Frequency 4.5e6 Hz} \\ 
\hline
 \multicolumn{1}{|p{1.5cm}|}{\centering \cellcolor{gray!30}F (mm)} & \multicolumn{1}{p{1.5cm}|}{\centering \cellcolor{orange!30}$M = 20$} & \multicolumn{1}{p{1.5cm}|}{\centering\cellcolor{yellow!20}$ M = 50$} & \multicolumn{1}{p{1.5cm}|}{\centering\cellcolor{orange!30}$M = 20$} & \multicolumn{1}{p{1.5cm}|}{\centering\cellcolor{yellow!20}$M = 50$} \\ 
\hline
10 & 5.83e-4 & 1.10e-4 & 9.03e-6 & 1.13e-5\\
\hline
25 & 1.37e-5 & 3.05e-5 & 1.04e-5 & 1.10e-5 \\
\hline
35 & 1.27e-3 & 1.15e-4 & 1.75e-5 & 7.28e-5\\
\hline
\end{tabular}
\caption{Quantitative measure $D$ for comparing $12$ experiments in the narrowband case. Each row contains the values for a different choice of frequency, number of active elements M, and depth of focus F.}
\label{tab:NarrowQuant}
\end{table}

\begin{figure}
    \centering
    \includegraphics[width=1.\textwidth]{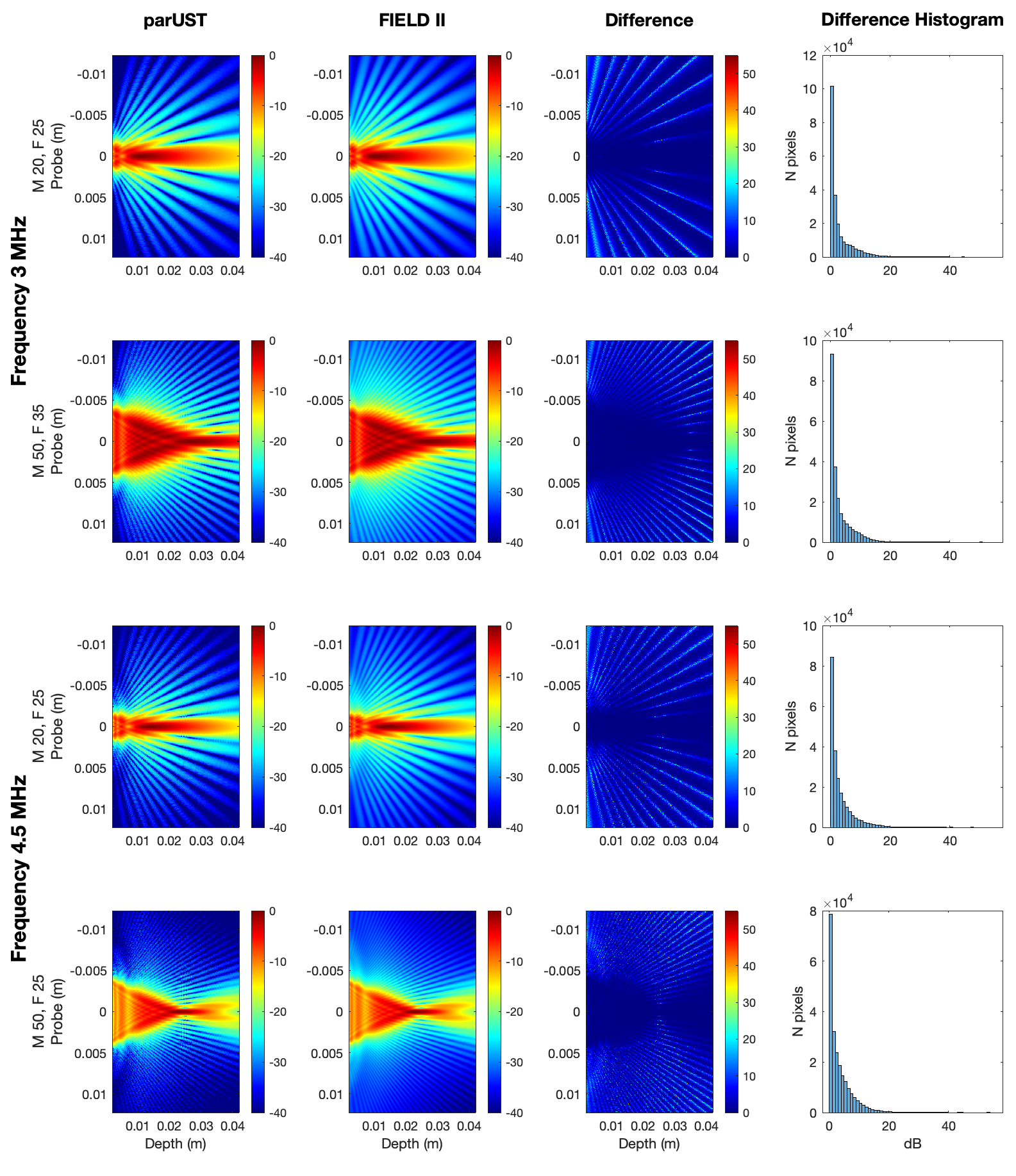}
    \caption{Comparison of generated narrowband BPs in different settings. On each column from the left: BPs generated with our simulator parUST, BPs from FIELD II, their absolute difference, and the distribution of pixel value differences. The unity of measure is dB.}
    \label{fig:NarrowPanel}
\end{figure}

These lower values correspond to very small differences between the beam patterns, as it is clear from Figure \ref{fig:NarrowPanel}.
The histogram of differerences in $dB$ presents some bins of relatively higher value with respect to the wideband case. However, this is due to the significant presence of very low intensity regions between the side lobes. In this case, even a very slight misalignment in side lobe position may produce a big difference in $dB$, that is however of negligible significance. 

\subsection{Evaluation of the software performance}
All tests have been performed on a laptop with $2.6\, GHz$ Intel Core i7 6 core processor and a $16\,Gb$ memory.
For a better performance on the computation of maps, many tests have been run on an nVidia workstation DGX-A100 with $128$ cores achieving the following average performance (Table \ref{tab:times}):
\begin{enumerate}
    \item $50 \, s$ for the wideband case, considering a maximum number of $2048$ time samples, thus obtaining $1024$ frequency samples saved for each point for the beam pattern computation on a grid of dimensions $[192 \cdot 4, 200]$;
    \item $16 \, s$ for the narrowband case, considering pre-computing and saving $35$ maps of dimension $[70 \cdot 4, 200]$ points, corresponding to $70$ probe elements. 
\end{enumerate}
The differences in computing times are consistent with the differences in the dimensions of the maps and with the fact that in the wideband case we need to compute the FFT of all impulse responses. 
It should be noted that the times needed for BP computation, dependent on variable transmit parameters, are significantly lower than the times needed for maps computation, which are conveniently performed only once given the probe model and the spatial grid. This difference becomes dramatic for the narrow band case. 
Moreover, if we compare the BP computation time with the global time of FIELD II we can clearly see the advantage of the proposed approach for optimization purposes, in which different sets of parameters need to be evaluated for many iterations with the same maps. 
It is worth noticing that the differences in computation time between the two devices for FIELD II are small due to the architecture characteristics. 
In detail, the laptop performs the computation over $6$ cores using the entire CPU, while the server, despite having many cores, uses just one if the code is not parallelized.
\begin{table}[!ht]
\centering 
\begin{tabular}{l|c c|c c| c c|}
\cline{2-6}
    & \multicolumn{2}{c|}{\cellcolor{red!30}Maps Computing}  & \multicolumn{2}{c|}{\cellcolor{red!30}BP Computing}& \multicolumn{2}{c|}{\cellcolor{red!30}FIELD II}\\ 
    \cline{2-6}  
    &\multicolumn{1}{p{1.3cm}|}{\cellcolor{orange!30}\centering Laptop} &\multicolumn{1}{p{1.3cm}|} {\cellcolor{yellow!30}\centering Server} &\multicolumn{1}{p{1.3cm}|}{\cellcolor{orange!30}\centering Laptop} & \multicolumn{1}{p{1.3cm}|}{\cellcolor{yellow!30}\centering Server} &\multicolumn{1}{p{1.3cm}|}{\cellcolor{orange!30}\centering Laptop} &
    \multicolumn{1}{p{1.3cm}|}{\cellcolor{yellow!30}\centering Server} \\ 
        \hline
\multicolumn{1}{|l|}{Wideband}   & \multicolumn{1}{c|}{720 s}       &     50 s   & \multicolumn{1}{c|}{8 s}  &   5 s  & \multicolumn{1}{c|}{80 s} &  76 s\\ \hline
\multicolumn{1}{|l|}{Narrowband} & \multicolumn{1}{c|}{246 s}       &   16 s & \multicolumn{1}{c|}{0.007 s}  & 0.004 s  & \multicolumn{1}{c|}{115 s}  &  105 s \\ \hline
\end{tabular}
\caption{Comparison of average computational times needed for the computation of maps and BPs, with both our simulator (first $4$ columns) and FIELD II (last $2$ columns). Once the maps are computed and stored, the running time for BP generation is significantly advantageous with respect to FIELD II.}
\label{tab:times}
\end{table}

\section{Conclusions}
The main contributions of this paper are: (1) a mathematical framework that deploys a hybrid time-frequency domain formulation for beam pattern modelling in the transmission step, based on which the parameters’ domain results in an N-dimensional torus, and (2) a novel simulator parUST \cite{parust} for beam pattern generation motivated by this framework, along with an experimental validation on simulated ultrasound data. 

The developed model has several advantageous characteristics. 
First, it is based on a two-step approach that allows one to pre-compute the system frequency response for a given probe model and region of interest: this methodology is particularly suited when the simulator is used as a kernel in gradient-based or non gradient-based optimization problems to the effect that many different sets of transmit parameters on the beam pattern shape have to be tested. 
Conversely, the reference simulator FIELD II has a monolithic structure that restarts all the computation every time a single parameter is changed. 
Second, our frequency domain approach avoids time oversampling, and it leads to time convolutions between impulse responses and waveform or attenuation terms being boiled down to element-wise products on a limited set of frequency bins.
Also, it specializes the narrowband case with a single frequency computation, hence gaining a substantial computational saving with respect to FIELD II-like time-domain approaches. The narrowband case formulation also allows one to investigate on the topology of space transmit delays, paving the way to suitable manifold optimization.
Finally, the simulator is an open-source Python package optimized with parallel computations, and it provides a high-performance trusted tool in view of a parameter optimization process on manifolds for BP formation and, ultimately, for PSF delineation.
It was numerically tested that the proposed model is both computationally efficient and compatible with realistic measurements and existing numerical simulators, and it may contribute to the design of new medical diagnostic ultrasound equipment.

Our investigations did not touch the problem of extending the model formulation to other types of probe geometries like convex or phased arrays, nor the possible extension of the to-be-optimized parameters to transmit apodization or element-dependent waveforms. Further, the receive and the image formation steps are yet to be considered. Such considerations are left for future studies.

\bibliographystyle{plain}
\bibliography{Razzetta}

\begin{thebibliography}{10}

\bibitem{Besson2019}
A.~Besson, J.~Thiran, and Y.~Wiaux.
\newblock Imaging from echoes: On inverse problems in ultrasound.
\newblock {\em PhD Thesis, Ecole Polytechnique F\'ed\'erale de Lausanne}, 2019.

\bibitem{Brody1974}
W.~R. Brody and J.~D. Meindl.
\newblock {Theoretical Analysis of the CW Doppler Ultrasonic Flowmeter}.
\newblock {\em IEEE Transactions on Biomedical Engineering},
  BME-21(3):183--192, 1974.

\bibitem{Bruyneel2013}
T.~Bruyneel, A.~Ortega, L.~Tong, and J.~D'hooge.
\newblock A {GPU}-based implementation of the spatial impulse response method
  for fast calculation of linear sound fields and pulse-echo responses of array
  transducers.
\newblock In {\em 2013 IEEE International Ultrasonics Symposium (IUS)}, pages
  367--369, 2013.

\bibitem{Chen2016}
Z.~Chen, A.~Basarab, and D.~Kouam\'e.
\newblock Compressive deconvolution in medical ultrasound imaging.
\newblock {\em IEEE Transactions on Medical Imaging}, 35(3):728–737, 2016.

\bibitem{Eldar2014}
T.~Chernyakova and Y.~C. Eldar.
\newblock Fourier-domain beamforming: the path to compressed ultrasound
  imaging.
\newblock {\em IEEE Transactions on Ultrasonics, Ferroelectrics, and Frequency
  Control}, 61(8):1252--1267, 2014.

\bibitem{Cobbold2007}
R.~S.~C. Cobbold.
\newblock {\em Foundations of Biomedical Ultrasound}.
\newblock Biomedical engineering series. Oxford University Press, 2007.

\bibitem{Crombie1997}
P.~Crombie, A.~Bascom, and S.~Cobbold.
\newblock {Calculating the pulsed response of linear arrays: Accuracy vs.
  computational efficiency}.
\newblock {\em IEEE Trans. Ultrason. Ferroelect. Freq. Contr.}, 44:997--1009,
  1997.

\bibitem{Cueto2022}
C.~Cueto, O.~Bates, G.~Strong, J.~Cudeiro, F.~Luporini, \`O.~Calder\'on Agudo,
  G.~Gorman, L.~Guasch, and M.~Tang.
\newblock Stride: A flexible software platform for high-performance ultrasound
  computed tomography.
\newblock {\em Computer Methods and Programs in Biomedicine}, 221:106855, 2022.

\bibitem{ARFI3}
J.~R. Doherty, G.~E. Trahey, K.~R. Nightingale, and M.~L. Palmeri.
\newblock Acoustic radiation force elasticity imaging in diagnostic ultrasound.
\newblock {\em IEEE transactions on ultrasonics, ferroelectrics, and frequency
  control}, 60(4):685--701, 2013.

\bibitem{Fenster2015}
A.~Fenster and J.~C. Lacefield.
\newblock {\em Ultrasound Imaging and Therapy}.
\newblock CRC Press, 1 2015.

\bibitem{Friis1937}
H.~T. Friis and C.~B. Feldman.
\newblock A multiple unit steerable antenna for short-wave reception.
\newblock {\em Proceedings of the Institute of Radio Engineers},
  25(7):841--917, 1937.

\bibitem{GARCIA2022}
D.~Garcia.
\newblock {SIMUS}: An open-source simulator for medical ultrasound imaging.
  part i: Theory \& examples.
\newblock {\em Computer Methods and Programs in Biomedicine}, 218:106726, 2022.

\bibitem{Basarab2022}
S.~Goudarzi, A.~Basarab, and H.~Rivaz.
\newblock A unifying approach to inverse problems of ultrasound beamforming and
  deconvolution.
\newblock {\em \url{https://arxiv.org/pdf/2112.14294.pdf}}, 2022.

\bibitem{Guo2022}
H.~Guo, H.~Xie, G.~Zhou, and N.~Q. Nguyen.
\newblock Fourier-domain beamforming and sub-nyquist sampling for coherent
  pixel-based ultrasound imaging.
\newblock In {\em 2022 IEEE International Ultrasonics Symposium (IUS)}, pages
  1--4, 2022.

\bibitem{Harris1981}
G.~R. Harris.
\newblock Transient field of a baffled planar piston having an arbitrary
  vibration amplitude distribution.
\newblock {\em The Journal of the Acoustical Society of America},
  70(1):186--204, 1981.

\bibitem{Hundt1980}
E.E. Hundt and E.A. Trautenberg.
\newblock Digital processing of ultrasonic data by deconvolution.
\newblock {\em IEEE Transactions on Sonics and Ultrasonics}, 27(5):249--252,
  1980.

\bibitem{Hyun2019}
D.~Hyun, L.~L. Brickson, K.~T. Looby, and J.~J. Dahl.
\newblock Beamforming and speckle reduction using neural networks.
\newblock {\em IEEE Transactions on Ultrasonics, Ferroelectrics, and Frequency
  Control}, 66(5):898--910, 2019.

\bibitem{Jensen1991}
J.~A. Jensen.
\newblock A model for the propagation and scattering of ultrasound in tissue.
\newblock {\em J Acoust Soc Am. 1991 Jan}, 89:182--191, 1991.

\bibitem{Jensen1996}
J.~A. Jensen.
\newblock {FIELD}: A program for simulating ultrasound systems.
\newblock {\em Medical and Biological Engineering and Computing}, 34:351--352,
  01 1996.

\bibitem{Jensen1999}
J.~A. Jensen.
\newblock {\em Linear description of ultrasound imaging systems: Notes for the
  International Summer School on Advanced Ultrasound Imaging at the Technical
  University of Denmark}.
\newblock Technical University of Denmark, Department of Electrical
  Engineering, 1999.

\bibitem{Jensen1999H}
J.~A. Jensen.
\newblock A new calculation procedure for spatial impulse responses in
  ultrasound.
\newblock {\em The Journal of the Acoustical Society of America},
  105(6):3266--3274, 06 1999.

\bibitem{Jensen2002}
J.~A. Jensen.
\newblock {\em Ultrasound Imaging and Its Modeling}, pages 135--166.
\newblock Springer Berlin Heidelberg, Berlin, Heidelberg, 2002.

\bibitem{Jensen2014}
J.~A. Jensen.
\newblock {A Multi-threaded Version of Field II}.
\newblock In {\em Proceedings of IEEE International Ultrasonics Symposium.
  IEEE}, pages 2229--2232, 2014.

\bibitem{Greenleaf1994}
J.~Y. Lu and J.~F. Greenleaf.
\newblock A study of two-dimensional array transducers for limited diffraction
  beams.
\newblock {\em IEEE Transactions on Ultrasonics, Ferroelectrics, and Frequency
  Control}, 41(5):724--39, 1994.

\bibitem{Luchies2018}
A.~C. Luchies and B.~C. Byram.
\newblock Deep neural networks for ultrasound beamforming.
\newblock {\em IEEE Transactions on Medical Imaging}, 37(9):2010--2021, 2018.

\bibitem{Michailovich2007}
O.~Michailovich and A.~Tannenbaum.
\newblock Blind deconvolution of medical ultrasound images: A parametric
  inverse filtering approach.
\newblock {\em IEEE Transactions on Image Processing}, 16(12):3005--3019, 2007.

\bibitem{Ozkan2018}
E.~Ozkan, V.~Vishnevsky, and O.~Goksel.
\newblock Inverse problem of ultrasound beamforming with sparsity constraints
  and regularization.
\newblock {\em IEEE Transactions on Ultrasonics, Ferroelectrics, and Frequency
  Control}, 65(3):356–365, 2018.

\bibitem{parust}
C.~Razzetta, M.~Crocco, and F.~Benvenuto.
\newblock {parallel parametric UltraSound Transmission software (parUST)}
  (version v1: 2023-05-29), {GitHub}:
  \url{https://github.com/chiararazzetta/parUST}, 2023.

\bibitem{Razzetta2023}
C.~Razzetta, M.~Crocco, and F.~Benvenuto.
\newblock A stochastic approach to delays optimization for narrowband transmit
  beam pattern in medical ultrasound, 2023.

\bibitem{Szabo2004}
T.~L. Szabo.
\newblock {\em Diagnostic ultrasound imaging: inside out}.
\newblock Academic press, 2004.

\bibitem{Basarab2016}
T.~Szasz, A.~Basarab, and D.~Kouam\'e.
\newblock Beamforming through regularized inverse problems in ultrasound
  medical imaging.
\newblock {\em IEEE Transactions on Ultrasonics, Ferroelectrics, and Frequency
  Control}, 63(12):2031--2044, 2016.

\bibitem{Tanter2014}
M.~Tanter and M.~Fink.
\newblock Ultrafast imaging in biomedical ultrasound.
\newblock {\em IEEE transactions on ultrasonics, ferroelectrics, and frequency
  control}, 61(1):102--119, 2014.

\bibitem{kwave2010}
B.~E. Treeby and B.~T. Cox.
\newblock {k-Wave: MATLAB toolbox for the simulation and reconstruction of
  photoacoustic wave fields}.
\newblock {\em Journal of Biomedical Optics}, 15(2):021314, 2010.

\bibitem{Villalba2019}
A.~I. Villalba, T.~Landry, and J.~Brown.
\newblock Parallel computing using {P}ython-based software for a high-frequency
  ultrasound system.
\newblock {\em The Journal of the Acoustical Society of America},
  146(4):3073--3073, 10 2019.

\bibitem{Li1999}
L.~Yadong and J.A. Zagzebski.
\newblock A frequency domain model for generating {B}-mode images with array
  transducers.
\newblock {\em IEEE Transactions on Ultrasonics, Ferroelectrics, and Frequency
  Control}, 46(3):690--699, 1999.

\bibitem{Yu2022}
B.~Yu, H.~Jin, Y.~Mei, J.~Chen, E.~Wu, and K.~Yang.
\newblock 3-{D} ultrasonic image reconstruction in frequency domain using a
  virtual transducer model.
\newblock {\em Ultrasonics}, 118:106573, 2022.

\end{thebibliography}
	
\end{document}